\documentclass[10pt,leqno]{article}
\usepackage{amsmath,amsthm}
\usepackage{amsfonts,amssymb,amsbsy}
\newtheorem{theorem}{Theorem}
\newtheorem{lemma}[theorem]{Lemma}
\newtheorem{proposition}[theorem]{Proposition}

\title{ Dyson's Brownian motions, intertwining and interlacing
}
\author{JON WARREN}
\date{}
\begin{document}
\maketitle
\begin{abstract}
A  family of reflected Brownian motions is used to construct Dyson's process of  non-colliding Brownian motions. A number of explicit 
formulae are given, including one for the distribution of a family of coalescing Brownian motions.
\end{abstract}
\section{Introduction}

The ordered eigenvalues $Y_1(t)\leq Y_2(t) \leq \ldots \leq Y_n(t)$  of a Brownian motion in the space of $ n \times n$ Hermitian matrices  
form a diffusion process which satisfies the stochastic differential equations,
\begin{equation}
\label{eqndyson}
Y_i(t)= y_i+ \beta_i(t)+\sum_{j \neq i}\int_0^t \frac{ds}{Y_i(s)-Y_j(s)},
\end{equation}
where $\beta_1,\beta_2,\ldots, \beta_n$ are independent real Brownian motions.
This is a  result that goes back to Dyson \cite{dyson} and we will refer to $Y$ as a Dyson non-colliding Brownian motion. 
A number of important papers in recent years have developed a  link between random matrices and certain combinatorial models, involving 
random permutations, last passage percolation, random tilings,  random growth models and queuing systems, see Baik, Deift and Johansson, 
\cite{baikdeiftjohansson} and Johansson,
\cite{johansson}, amongst many others. A recent survey is given by K\"{o}nig \cite{konig}. At the heart of this connection lies  the 
Robinson-Schensted-Knuth algorithm, a combinatorial procedure which has its origins  in group representation theory, and using this the 
following remarkable formula, was observed by Gravner, Tracy and Widom, \cite{gravnertracywidom} and Baryshnikov \cite{baryshnikov}, 
representing the largest eigenvalue $Y_n(t)$ ( assuming $Y(0)=0$) in terms of independent, real-valued, Brownian motions $B_1,B_2, \ldots, 
B_n$, 
\begin{equation}
\label{identity}
Y_n(t) \stackrel{dist}{=} \sup_{0=t_0 \leq t_1 \leq  \ldots \leq t_n=t} \sum_{i=1}^n \bigl\{ B_i(t_i)-B_{i}(t_{i-1})\bigr\}.
\end{equation}
O'Connell and Yor, \cite{oconnellyor}, give a proof of this identity by considering reversibility properties of a queuing system, which in a 
subsequent paper, O'Connell \cite{oconnell}, is shown to be linked to the RSK algorithm also. 
Another proof, again involving RSK, is given 
by Doumerc, \cite{doumerc}.

In this paper an  a different proof of  the identity \eqref{identity} is given,  based around the following construction.  
Let $\bigl(Y(t);t \geq 0\bigr)$ be a Dyson process, with components $Y_1, Y_2, \ldots Y_n$ solving \eqref{eqndyson}. Let $\bigl(X(t);t \geq 
0 \bigr)$ be  a process with $(n+1)$ components which  are {\em interlaced } with those of  $Y$, meaning 
\begin{equation}
X_1(t) \leq Y_1(t)\leq X_2(t) \ldots \leq Y_n(t) \leq X_{n+1}(t), \qquad \text{ for all } t \geq 0,
\end{equation}
and which satisfies the equations
\begin{equation}
\label{eqnref}
X_i(t)= x_i+\gamma_i(t) +  \bigl\{L^-_i(t)- L^+_i(t)\bigr\}.
\end{equation}
Here $\bigl(\gamma(t); t\geq 0\bigr)$ is a standard Brownian motion in ${\mathbf R}^{n+1}$, independent of the Brownian motion $\beta$ which 
drives $Y$. The processes $(L_i^+(t); t \geq 0)$ and $(L_i^-(t);t \geq 0 )$ are continuous non-decreasing processes that increase only at 
times when $X_i(t)=Y_{i}(t)$ and $X_i(t)=Y_{i-1}(t)$ respectively: they are twice the semimartingale local times at zero of $X_i-Y_i$ and 
$X_i-Y_{i-1}$. The two exceptional cases $L^-_1(t)$ and $L^+_{n+1}(t)$ are defined to be identically zero.
Conditionally on $Y$ the particles corresponding to $X$ evolve as independent 
Brownian motions except when collisions occur with particles corresponding to  $Y$. Think of the particles corresponding to the components 
of $Y$ as being ``heavy'' so that in collisions with the ``light'' particles corresponding to components of $X$  their motion is unaffected. 
On the other hand the light particles receive a singular drift from the collisions which maintains the interlacing.
We will verify that is possible to start $X$ and $Y$ from the origin so that $x_i=y_j=0$ for all $i$ and $j$.  Then, see Proposition 
\ref{interlace}, {\em the process $X$ is distributed as a  Dyson non-colliding process with $(n+1)$ particles.}
Thus if we  observe only the  particles corresponding to the components of $X$, the singular drifts that these particles experience from 
collisions with the unseen particles corresponding to $Y$ are somewhat magically  transmuted into an electrostatic repulsion. This is a 
consequence of a relationship between the semigroup of the extended process $(X,Y)$ and the semigroup of $X$ that is called an {\em 
intertwining } relation.

\section{A duality between interlaced Brownian motions}

Consider a  continuous, adapted, ${\mathbf R}^{n+1}\times {\mathbf R}^n$-valued process $\bigl(X(t),Y(t);t \geq 0\bigr)$ having components 
$X_1(t),X_2(t), \ldots X_{n+1}(t)$ and $Y_1(t),Y_2(t),\ldots, Y_n(t)$ which is defined on filtered probability space $\bigl(\Omega,{\cal 
F},\bigl\{{\cal F}_t\bigr\}_{t \geq 0}, {\mathbf Q}^n_{x,y}\bigr)$  satisfying, for all $t \geq 0$,  the interlacing condition
\[
X_1(t) \leq Y_1(t) \leq X_2(t)\leq \ldots \leq Y_n(t) \leq X_{n+1}(t),
\]
and the equations
\begin{align}
\label{eqnfory}
Y_i(t)&=y_i+\beta_i(t\wedge \tau), \\
\label{eqnforx}
X_i(t)&=x_i+\gamma_i(t \wedge \tau)+  L_i^-(t\wedge \tau)-L_i^+(t\wedge \tau),
\end{align}
where,
\begin{description}
\item $\tau$ is the stopping time given by $\tau=\inf\bigl\{ t\geq 0: Y_i(t)=Y_{i+1}(t) \text{ for some } i\in \{1,2,\ldots,n-1\}\bigr\}$,
\item $ \beta_1,\beta_2, \ldots \beta_{n}, \gamma_1,\gamma_2, \ldots \gamma_{n+1}$ are independent ${\cal F}_t$-Brownian motions,
\item $L_1^-(t)=L_{n+1}^+(t)=0$ for all $t \geq 0$, otherwise the processes $L^+_i$ and $L^-_i$  are continuous, non-decreasing and increase 
only when $X_i=Y_{i}$ and $X_i=Y_{i-1}$ respectively, 
\[
L^+_i(t)= \int_0^t {\mathbf 1}\bigl(X_i(s)=Y_i(s)\bigr) dL^+_i(s) \qquad L^-_i(t)= \int_0^t {\mathbf 1}\bigl(X_i(s)=Y_{i-1}(s)\bigr) 
dL^-_i(s).
\]
\end{description}
The process just defined is called a stopped, semimartingale reflecting Brownian motion. For general results on such processes see, for 
example, Dai and Williams, \cite{daiwilliams}. In this case it is not difficult to give a pathwise construction starting from the Brownian 
motions $\beta_i$, for $i \in \{1,2,\ldots, n\}$, and $\gamma_i$ for $i \in \{1,2,\ldots, n+1\}$, together with the choice of initial 
co-ordinates  $x_1 \leq y_1 \leq x_2 \ldots \leq y_n \leq x_{n+1}$.  We obtain $Y_i$ immediately. $X_i$ is constructed by alternately using 
the usual Skorokhod construction to push $X_i$ up from $Y_{i-1}$ and down from $Y_{i}$. For more details  see Section 3 of  
\cite{soucaliuctothwerner}, where  a similar construction is used. In fact by the same argument as Lemma 6 in Soucaliuc, Toth and Werner, 
\cite{soucaliuctothwerner}  pathwise uniqueness holds, and hence the law of $\bigl(X,Y\bigr)$ is uniquely determined.
This uniqueness implies, by standard methods, that the process is Markovian, and in fact we are able to give an explicit formula for its 
transition probabilities. 

 We denote by $\phi_t$ the centered Gaussian density with variance $t$. $\Phi_t$ is the corresponding distribution function 
\[
\Phi_t(y)= \int_{-\infty}^y \frac{1}{\sqrt{2\pi t}} \exp\{-z^2/(2t)\}dz,
\]
and 
\[
\phi^\prime_t(y)= \frac{ -y}{\sqrt{2\pi t^3}} \exp\{-y^2/(2t)\}.
\]
Let $W^{n+1,n}=\{ (x,y) \in {\mathbf R}^{n+1} \times {\mathbf R}^n: x_1 \leq y_1 \leq x_2\leq \ldots \leq y_n \leq x_{n+1} \}$. 
Define $q^n_t\bigl((x,y), (x^\prime, y^\prime) \bigr)$ for $(x,y), (x^\prime,y^\prime) \in W^{n+1,n}$ and $t>0$ to be equal to determinant 
of the $(2n+1)\times(2n+1)$ matrix 
\[
 \begin{pmatrix}
A_t( x,x^\prime) & B_t(x, y^\prime) \\
C_t(y,x^\prime)  & D_t(y,y^\prime)
\end{pmatrix}
\]
where
\begin{description}
\item  $A_t(x,x^\prime)$ is an $(n+1) \times (n+1)$ matrix with $(i,j)$th element $\phi_t(x^\prime_j-x_i)$;
\item   $B_t(x,y^\prime)$ is an $(n+1) \times n $ matrix with $(i,j)$th element $\Phi_t(y^\prime_j-x_i)-{\mathbf 1}(j\geq i)$. 
\item   $C_t(y,x^\prime)$ is an $n \times (n+1) $ matrix with $(i,j)$th element $\phi^\prime_t(x^\prime_j-y_i)$;
\item   $D_t(y,y^\prime)$ is an $n \times n $ matrix with $(i,j)$th element $\phi_t(y^\prime_j-y_i)$.
\end{description}

\begin{lemma}
\label{weakcon}
 For any $f:W^{n+1,n} \rightarrow {\mathbf R}$ which is bounded and continuous, and zero in a neighbourhood of the boundary of $W^{n+1,n}$,
\[
\lim_{t \downarrow 0} \int_{W^{n+1,n}} q^n_t\bigl( w,w^\prime\bigr)f(w^\prime)dw^\prime = f(w),
\]
uniformly  for all $w=(x,y) \in W^{n+1,n}$.
\end{lemma}
\begin{proposition} 
\label{trans1}
$(q^n_t; t >0)$ are a family of transition densities for   the  process $\bigl(X,Y\bigr)$ killed at the instant $\tau$, that is to say,  for 
$t>0$ and  $(x,y),(x^\prime,y^\prime)\in W^{n+1,n}$,
\[q^n_t\bigl((x,y),(x^\prime,y^\prime)\bigr)dx^\prime dy^\prime={\mathbf Q}^n_{x,y} \bigl( X_t \in dx^\prime, Y_t\in dy^\prime; t < \tau).
\]
\end{proposition}

\begin{proof}
For any choice of $z^\prime \in {\mathbf R}$, each of the functions $ (t,z) \mapsto \Phi_t (z^\prime-z)$, $(t,z)\mapsto \phi_t(z^\prime-z)$ 
and $(t,z)\mapsto \phi_t^\prime(z^\prime-z)$ satisfies the heat equation on $(0,\infty)\times {\mathbf R}$. Thus, by  differentiating  the 
determinant, we find that,
\begin{equation}
\tfrac{1}{2}  \sum_{i=1}^{2n+1} \frac{\partial^2q^n_t}{\partial w_i^2} (w, w^\prime) = \frac{\partial q^n_t}{\partial t}(w, w^\prime)   
\quad  (t,w,w^\prime)\in (0,\infty)\times {\mathbf R}^{2n+1}\times {\mathbf R}^{2n+1}.
\end{equation}
We need to identify certain boundary conditions.  We treat $w^\prime=(x^\prime,y^\prime)\in W^{n+1,n}$ as fixed.   First consider $(x,y) \in 
\partial W^{n+1,n}$ satisfying  $y_i=y_{i+1}$ for some $i\in \{1,2,\ldots n-1\}$.  We see that the $i$th and $(i+1)$th rows of both 
$C_t(y,x^\prime)$ and $D_t(y,y^\prime)$ are equal, and hence 
$q^n_t\bigl((x,y),(x^\prime,y^\prime)\bigr)=0$. 
 Next consider $(x,y)\in \partial W^{n+1,n}$ satisfying $x_i= y_i$ for some $i \in \{1,2,\ldots n\}$.  Calculate $
\frac{\partial}{\partial x_i} q^n_t( (x,y), (x^\prime,y^\prime))$
by differentiating the $i$th rows of $A_t(x,x^\prime)$ and $B_t(x,y^\prime)$. Notice that, under our assumption that $x_i=y_i$, the $i$th 
row of $\frac{\partial}{\partial x_i} A_t(x,x^\prime)$ is equal to the $i$th row of $-C_t(y,x^\prime)$. Likewise   the $i$th row of 
$\frac{\partial}{\partial x_i} B_t(x,y^\prime)$ is equal to the $i$th row of $-D_t(y,x^\prime)$. 
Thus we deduce that $\frac{\partial}{\partial x_i} q^n_t( (x,y), (x^\prime,y^\prime))=0$. Finally consider  $(x,y) \in \partial W^{n+1,n}$ 
satisfying $x_{i+1}=y_i$ for some $i \in \{1,2,\ldots n\}$. Similarly to the previous case we obtain
$\frac{\partial}{\partial x_{i+1}} q^n_t( (x,y), (x^\prime,y^\prime))=0$.

Let $f: W^{n+1,n} \rightarrow {\mathbf R}$ be a bounded and continuous, and zero in a neighbourhood of the boundary.  Then   define a smooth 
function $F$  on $(0,\infty)\times W^{n+1,n}$ via
\begin{equation*}
F(t,w)= \int_{W^{n+1,n}} q^n_t (w, w^\prime) f(w^\prime)dw^\prime.
\end{equation*}
By virtue of the above observations regarding $q^n_t$, and differentiating through the integral,  we find that
\begin{equation*}
\tfrac{1}{2} \sum_{i=1}^{2n+1} \frac{\partial^2 F}{\partial w_i^2} (t,w)
 = \frac{\partial F}{\partial t}(t, w)   \quad \text{ on } (0,\infty)\times W^{n+1,n},
\end{equation*}
with the boundary conditions
\begin{align*}  &F(t,w)=0 \text{  whenever $ w=(x,y)$ satisfies 
$y_i=y_{i+1}$} \\
&\frac{\partial F}{\partial x_i}(t,w)=0 \text{  whenever $ w=(x,y)$ satisfies 
$x_i=y_{i}$} \\
&\frac{\partial F}{\partial x_{i+1}}(t,w)=0 \text{  whenever $ w=(x,y)$ satisfies 
$x_{i+1}=y_{i}$} 
\end{align*}
Fix $T, \epsilon>0$. Applying It\^{o}'s formula, we find that the process $\bigl(F\bigl((T+\epsilon-t, (X_t,Y_t)\bigr); t\in [0,T]\bigr)$ is 
a local martingale, which is easily seen to be bounded and hence is a true martingale. Thus 
\begin{align*}
F\bigl(T+\epsilon, (x,y) \bigr)&={\mathbf Q}^n_{x,y} \bigl[ F\bigl(\epsilon,  (X_{T},Y_T)\bigr) \bigr]\\
&={\mathbf Q}^n_{x,y} \bigl[ F\bigl(\epsilon,  (X_{T},Y_T)\bigr){\mathbf 1}(T<\tau) \bigr].
\end{align*}
Appealing to the previous lemma, we may let $\epsilon \downarrow 0$ and so obtain,
\[
F\bigl(T, (x,y)\bigr)= {\mathbf Q}^n_{x,y} \bigl[ f(X_{T},Y_T)){\mathbf 1}(T<\tau) \bigr].
\]
Since the part of the distribution of $(X_T, Y_T)$ that charges the boundary of $W^{n+1,n}$ exactly corresponds to the event $\{T \geq 
\tau\}$ this
suffices to prove the proposition.
\end{proof}

We now consider a second reflected semimartingale Brownian motion$\bigl(\hat{X},\hat{Y}\bigr)$  having components 
$\hat{X}_1(t),\hat{X}_2(t), \ldots \hat{X}_{n+1}(t)$ and $\hat{Y}_1(t),\hat{Y}_2(t),\ldots, \hat{Y}_n(t)$ which is defined on filtered 
probability space $\bigl(\Omega,{\cal F},\bigl\{{\cal F}_t\bigr\}_{t \geq 0}, \hat{{\mathbf Q}}^n_{x,y}\bigr)$  satisfying, for all $t \geq 
0$,  the interlacing condition
\[
\hat{X}_1(t) \leq \hat{Y}_1(t) \leq \hat{X}_2(t)\leq \ldots \leq \hat{Y}_n(t) \leq \hat{X}_{n+1}(t),
\]
and the equations
\begin{align}
\hat{Y}_i(t)&=y_i+\beta_i(t\wedge \hat{\tau})+L_i^-(t\wedge \hat{\tau})-L_i^+(t\wedge \hat{\tau}), \\
\hat{X}_i(t)&=x_i+\gamma_i(t \wedge \hat{\tau}),
\end{align}
where,
\begin{description}
\item $\hat{\tau}$ is the stopping time given by $\hat{\tau}=\inf\bigl\{ t\geq 0: \hat{X}_i(t)=\hat{X}_{i+1}(t) \text{ for some } i\in 
\{1,2,\ldots,n\}\bigr\}$,
\item $ \beta_1,\beta_2, \ldots \beta_{n}, \gamma_1,\gamma_2, \ldots \gamma_{n+1}$ are independent ${\cal F}_t$-Brownian motions,
\item  the processes $L^+_i$ and $L^-_i$  are continuous, non-decreasing and increase only when $\hat{Y}_i=\hat{X}_{i+1}$ and 
$\hat{Y}_i=\hat{X}_{i}$ respectively, 
\[
L^+_i(t)= \int_0^t {\mathbf 1}\bigl(\hat{Y}_i(s)=\hat{X}_{i+1}(s)\bigr) dL^+_i(s) \qquad L^-_i(t)= \int_0^t {\mathbf 
1}\bigl(\hat{Y}_i(s)=\hat{X}_{i}(s)\bigr) dL^-_i(s).
\]
\end{description}
Notice the difference between this process and $\bigl(X,Y\bigr)$ is the reflection rule:  here $\hat{Y}$ is pushed off $\hat{X}$ whereas it 
was $X$ that was pushed off $Y$.

Define a family $\bigl(\hat{q}^n_t; t > 0\bigr)$ via
\begin{equation}
\label{dual}
\hat{q}^n_t\bigl((x,y),(x^\prime,y^\prime)\bigr)=q^n_t\bigl((x^\prime,y^\prime),(x,y)\bigr) \quad \text{ for } (x,y),(x^\prime,y^\prime)\in 
W^{n+1,n}.
\end{equation}
The following proposition is proved by arguments exactly parallel to those  just given in proof of Proposition \ref{trans1}.

\begin{proposition} $(\hat{q}^n_t; t > 0)$ are a family of transition densities for   the  process $\bigl(\hat{X},\hat{Y}\bigr)$ killed at 
the instant $\hat{\tau}$, that is to say,  for $t>0$ and  $(x,y),(x^\prime,y^\prime)\in W^{n+1,n}$,
\[\hat{q}^n_t\bigl((x,y),(x^\prime,y^\prime)\bigr)dx^\prime dy^\prime=\hat{\mathbf Q}^n_{x,y} \bigl( \hat{X}_t \in dx^\prime, \hat{Y}_t\in 
dy^\prime; t < \hat{\tau}).
\]
\end{proposition}

The duality, represented by \eqref{dual}, between the transition semigroups of $\bigl(X,Y\bigr)$ and $\bigl(\hat{X},\hat{Y}\bigr)$ is not 
unexpected.
It is consistent with general results, see for example DeBlassie \cite{deblassie}, and Harrison and Williams \cite{harrisonwilliams}, which 
show that, in a variety of contexts, the dual of a reflected Brownian motion is another reflected Brownian motion where the direction of 
reflection at the boundary is obtained by reflecting the original direction of reflection  across the normal vector. This is precisely the 
relationship holding between $\bigl(X,Y\bigr)$ and $\bigl(\hat{X},\hat{Y}\bigr)$  here.

\section{An intertwining involving Dyson's Brownian motions}

It is  known that Dyson's non-colliding Brownian motions can be obtained by means of a Doob $h$-transform. Let $W^n=\{y \in {\mathbf R}^n: 
y_1 \leq y_2 \leq \ldots \leq y_n \}$. Suppose that $\bigl(Y_t;t \geq 0 \bigr)$ is when governed by the probability measure ${\mathbf 
P}^n_y$ a standard  Brownian motion in ${\mathbf R}^n$, relative to a filtration $\{{\cal F}_t; t \geq 0\}$, starting  from  a point $y \in 
W^n$ and stopped at the instant $\tau=\inf\{t \geq 0: Y_i(t)=Y_j(t) \text{ for some } i \neq j \}$. The transition probabilities of $Y$ 
killed at the time $\tau$ are given explicitly by the Karlin-McGregor formula, \cite{karlinmcgregor},
\begin{equation}
{\mathbf P}^n_y\bigl( Y_t\in dy^\prime; t< \tau \bigr) = p^n_t(y,y^\prime) dy^\prime  ,
\end{equation}
for $y, y^\prime \in W^n$, where, with $\phi_t$ again denoting the Gaussian kernel with variance $t$,
\begin{equation}
p^n_t(y,y^\prime)= \det \bigl\{ \phi_t(y^\prime_j-y_i); 1 \leq i,j \leq n \bigr\}.
\end{equation}
If the initial co-ordinates $y$ satisfy $y_1<y_2< \ldots < y_n$, then we may define a new probability measure  by the absolute continuity 
relation
\begin{equation}
{\mathbf P}^{n,+}_{y} = \frac{h_n(Y_{t\wedge \tau})}{h_n(y)} \cdot {\mathbf P}^n_{y} \quad \text{ on } {\cal F}_t,
\end{equation}
for $t>0$, where $h_n$ is the function given by
\begin{equation}
h_n(y)= \prod_{i<j} (y_j-y_i).
\end{equation}
Under ${\mathbf P}^{n,+}_y$  the process $Y$ evolves as a Dyson non-colliding Brownian motion, that is to say $\tau$ is almost surely 
infinite and  the stochastic differential equations \eqref{eqndyson} hold. The transition probabilities
\begin{equation}
{\mathbf P}^{n,+}_y\bigl( Y_t\in dy^\prime; t< \tau \bigr) = p^{n,+}_t(y,y^\prime) dy^\prime,
\end{equation}
are related to those for the killed process by an $h$-transform
\begin{equation}
p_t^{n,+}(y,y^\prime)= \frac{h_n(y^\prime)}{h_n(y)}p^n_t(y,y^\prime),  
\end{equation}
for $y,y^\prime \in W^n\setminus \partial W^n$. 
Finally we recall, see O'Connell and Yor, \cite{oconnellyor}, that we  may describe  ${\mathbf P}_0^{n,+}$, the measure under which the 
non-colliding Brownian motion issues from the origin by  specifying  that it is Markovian with transition densities $\bigl( p_t^{n,+}; t>0 
\bigr)$ and with the entrance law 
\begin{equation}
{\mathbf P}_0^{n,+} \bigl( Y_t \in dy\bigr) = \mu^n_t(y)dy, 
\end{equation}
for $t>0$,  given by
\begin{equation}
\mu^n_t(y)= \frac{1}{Z_n} t^{-n^2/2} \exp \left\{ -\sum_{i} y_i^2/(2t)\right\} \left\{ \prod_{i<j} (y_j-y_i) \right\}^2,
\end{equation}
with the normalizing constant being $Z_n= (2\pi)^{n/2}\prod_{j<n} j!$.

Now suppose that $\bigl(X,Y\bigr)$ is governed by the probability measure ${\mathbf Q}^n_{x,y}$ defined in the previous section.   Recall 
that $\bigl(q^n_t; t > 0 \bigr)$ are the transition densities of the process killed at the time $\tau=\inf\{t \geq 0: Y_i(t)=Y_j(t) \text{ 
for some } i \neq j \}$. Suppose the initial co-ordinates $y$ of $Y$ satisfy $y_1<y_2< \ldots < y_n$, then we may define a new probability 
measure ${\mathbf Q}_{x,y}^{n,+}$ by the absolute continuity relation
\begin{equation}
{\mathbf Q}^{n,+}_{x,y} = \frac{h_n(Y_{t\wedge \tau})}{h_n(y)} \cdot {\mathbf Q}^n_{x,y} \quad \text{ on } {\cal F}_t,
\end{equation}
for $t>0$.
It follows from  the fact that under ${\mathbf Q}^n_{x,y}$ the process $Y$ evolves as a Brownian motion stopped at the instant $\tau$,  that 
$h_n(Y_{t\wedge \tau})$ is a martingale, and that this definition is hence consistent as $t$ varies. 
Under the measure ${\mathbf Q}^{n,+}_{x,y}$, the process $Y$  now evolves as a non-colliding Brownian motion satisfying the stochastic 
differential equation \eqref{eqndyson}, whilst the process $X$ satisfies \eqref{eqnref}.
 The corresponding transition densities $\bigl(q^{n,+}_t;t >0 \bigr)$ are obtained from those for the killed process  by the $h$-transform
\begin{equation}
q_t^{n,+}\bigl((x,y),(x^\prime,y^\prime)\bigr)= \frac{h_n(y^\prime)}{h_n(y)} q^n_t\bigl((x,y),(x^\prime,y^\prime)\bigr)
\end{equation}
for $(x,y),(x^\prime,y^\prime) \in W^{n+1,n}$ with the components of $y$ all distinct.
\begin{lemma}
\label{entrance}
The family of  probability measures with densities given by $\bigl(\nu^n_t; t>0\bigr)$ on $W^{n+1,n}$, given by
\begin{equation*}
\nu^n_t(x,y)=  \frac{n!}{Z_{n+1}} t^{-(n+1)^2/2} \exp \left\{ -\sum_{i} x_i^2/(2t)\right\} \left\{ \prod_{i<j} (x_j-x_i) \right\}\left\{ 
\prod_{i<j} (y_j-y_i) \right\},
\end{equation*}
form an entrance law for $\bigl(q^{n,+}_t; t>0\bigr)$, that is to say, for $s,t>0$
\[
\nu^n_{t+s}(w^\prime)= \int_{W^{n+1,n}} \nu^n_s(w)q_t^{n,+}(w,w^\prime)dw.
\]
\end{lemma}
Accordingly we may define a probability measure ${\mathbf Q}^{n,+}_{0,0}$, under which the process $\bigl(X,Y\bigr)$ is Markovian
 with transition densities  $\bigl( q_t^{n,+}; t>0 \bigr)$ and with the entrance law 
\begin{equation}
{\mathbf Q}_{0,0}^{n,+} \bigl( X_t\in dx ,Y_t \in dy \bigr) = \nu^n_t(x,y)dxdy.
\end{equation}
It is easy to see that under this measure $\bigl(X,Y\bigr)$ satisfies the equations \eqref{eqndyson} and \eqref{eqnref}, starting  from the 
origin $x=0, y=0$. Presumably any solution to \eqref{eqndyson} and \eqref{eqnref} starting from the origin  has the same law, but we do not 
prove this.

We may now state the main result of this section.

\begin{proposition}
\label{interlace}
Suppose the process $\bigl(X_t,Y_t; t \geq 0 \bigr)$ is governed by ${\mathbf Q}_{0,0}^{n,+}$ then the process $\bigl(X_t;t \geq 0\bigr)$ is 
distributed as
under ${\mathbf P}_0^{n+1,+}$, that is as a Dyson non-colliding Brownian motion  in $W^{n+1}$  starting from the origin.
\end{proposition}

This  result is proved by means of a criterion described by Rogers and Pitman \cite{rogerspitman} for  a function of a Markov process to be 
Markovian, see Carmona, Petit and Yor, \cite{carmonapetityor}, for futher examples of intertwinings.  For $x \in W^{n+1}$ let $W^n(x)=\{y 
\in {\mathbf R}^{n}: x_1 \leq y_1\leq \ldots \leq y_n\leq x_{n+1}$, and  define 
\begin{equation}
\lambda^n(x, y)= n!\frac{h_n(y)}{h_{n+1}(x)}, 
\end{equation}
for $x \in W^{n+1}\setminus \partial W^{n+1}$ and $ y \in W^{n}(x)$. 
The normalizing constant being chosen so that $\lambda^n( x, \cdot)$ is the density of a probability measure on $W^n(x)$. This follows from 
the equality
\begin{equation}
 \int_{W^n(x)} h_n(y)dy= \frac{1}{n!}h_{n+1}(x),
\end{equation}
which is easily verified by writing $h_n(y)=\det \bigl\{ y_i^{j-1}; 1\leq i,j \leq n \bigr\}$. 
The proof of Proposition \ref{interlace} depends on the  following intertwining relation between  $\bigl(q_t^{n,+}; t>0 \bigr)$ and 
$\bigl(p_t^{n+1,+}; t >0 \bigr)$,  for all $t>0$, $ x \in W^{n+1}\setminus \partial W^{n+1}$, and $(x^\prime,y^\prime)\in W^{n+1,n}$,
\begin{equation}
\label{intertwine}
\int_{W^n(x)} \lambda^n(x,y) q^{n,+}_t\bigl((x,y), (x^\prime,y^\prime)\bigr)dy= p_t^{n+1,+}(x,x^\prime)\lambda^n(x^\prime,y^\prime).
\end{equation}
This by be verified directly using the explicit formula for  $q^n_t$ given in the previous section. Alternatively the following derivation 
is enlightening.  
Recall that if $\bigl(\hat{X}_t,\hat{Y}_t; t \geq 0 \bigr)$ is governed by $\hat{\mathbf Q}^n_{x,y}$ then the process $\bigl(\hat{X}_t; 
t\geq 0 \bigr)$ is a Brownian motion stopped at the instant $\hat{\tau}=\inf\bigl\{t \geq 0; \hat{X}_i=\hat{X}_j \text{ for some } i \neq j 
\bigr\}$. Consequently the transition probabilities of the killed process satisfy
\begin{equation}
\int_{W^n(x^\prime)} \hat{q}^n_t\bigl( (x,y), (x^\prime,y^\prime)\bigr) dy^\prime =p^{n+1}_t(x,x^\prime).
\end{equation} 
Now using the duality between $q^n_t$ and $\hat{q}^n_t$ and the symmetry of $p^{n+1}_t$ we may re-write this as
\begin{equation}
\int_{W^n(x^\prime)} {q}^n_t\bigl( (x^\prime,y^\prime), (x,y)\bigr) dy^\prime =p^{n+1}_t(x^\prime,x).
\end{equation}
Finally to  obtain  \eqref{intertwine} we swop  the roles of $(x,y)$ and $(x^\prime,y^\prime)$ and use the expressions for $q_t^{n,+}$ and  
$p_t^{n+1,+}$ as $h$-transforms.
As a first application of the intertwining we have the following.
\begin{proof}[Proof of Lemma \ref{entrance}]
Notice that $\nu^n_t(x,y)= \mu^{n+1}_t(x)\lambda^n(x,y)$. Hence, by virtue of the intertwining and the fact that $\bigl(\mu^{n+1}_t; t 
>0\bigr)$ is an entrance law for $\bigl(p^{n+1,+}_t; t >0\bigr)$ we have,
\begin{multline*}
\int_{W^{n+1,n}} dxdy\; \nu^n_s(x,y)q^{n,+}_t\bigl((x,y),(x^\prime,y^\prime)\bigr)= \\
\int_{W^{n+1}} dx \;\mu_s^{n+1}(x) \int_{W^{n}(x)} dy\; \lambda^n(x,y)q^{n,+}_t\bigl((x,y),(x^\prime,y^\prime)\bigr) = \\
\int_{W^{n+1}} dx \;\mu_s^{n+1}(x) p^{n+1,+}_t(x,x^\prime)\lambda^n(x^\prime,y^\prime) =
\mu_{t+s}^{n+1}(x^\prime)\lambda^n(x^\prime,y^\prime)=\nu^n_{t+s}(x^\prime,y^\prime).
\end{multline*}

\end{proof}

A similar argument, following \cite{rogerspitman} proves the proposition.
\begin{proof}[Proof of Proposition \ref{interlace}]
For a sequence of times $0<t_1<t_2<\ldots <t_n$, repeated use of the intertwining relation gives,
\begin{multline*}
{\mathbf Q}_{0,0}^{n,+}\bigl(X_{t_1}\in A_1, X_{t_2}\in A_2, \ldots ,X_{t_n}\in A_n\bigr) =\\
\shoveleft{\int_{A_1} dx_1\ldots \int_{A_n} dx_n\int_{W^n(x_1)} dy_1\ldots \int_{W^n(x_n)}dy_n\;\nu^n_{t_1}(x_1,y_1) 
q^{n,+}_{t_2-t_1}\bigl((x_1,y_1),(x_2,y_2)\bigr)\ldots}\\
\shoveright{  \ldots q^{n,+}_{t_n-t_{n-1}}\bigl((x_{n-1},y_{n-1}),(x_n,y_n)\bigr)=}\\
 \shoveleft{\int_{A_1} dx_1\ldots \int_{A_n} dx_n\int_{W^n(x_1)} dy_1\ldots \int_{W^n(x_n)}dy_n\;\mu^{n+1}_{t_1}(x_1)\lambda^n(x_1,y_1) 
q^{n,+}_{t_2-t_1}\bigl((x_1,y_1),(x_2,y_2)\bigr)\ldots}\\
 \shoveright{ \ldots q^{n,+}_{t_n-t_{n-1}}\bigl((x_{n-1},y_{n-1}),(x_n,y_n)\bigr)=}\\
\shoveleft{ \int_{A_1} dx_1\ldots \int_{A_n} dx_n\int_{W^n(x_2)} dy_2\ldots \int_{W^n(x_n)}dy_n\;
 \mu^{n+1}_{t_1}(x_1) p^{n+1,+}_{t_2-t_1}\bigl(x_1,x_2\bigr)\lambda^n(x_2,y_2) \ldots} \\
 \shoveright{\ldots q^{n,+}_{t_n-t_{n-1}}\bigl((x_{n-1},y_{n-1}),(x_n,y_n)\bigr)=}\\
 \shoveleft{\int_{A_1} dx_1\ldots \int_{A_n} dx_n \int_{W^n(x_n)}dy_n\;
 \mu^{n+1}_{t_1}(x_1) p^{n+1,+}_{t_2-t_1}\bigl(x_1,x_2\bigr) \ldots p^{n+1,+}_{t_n-t_{n-1}}\bigl(x_{n-1},x_n\bigr)\lambda^n(x_n,y_n)=}\\
 \int_{A_1} dx_1\ldots \int_{A_n} dx_n \;
 \mu^{n+1}_{t_1}(x_1) p^{n+1,+}_{t_2-t_1}\bigl(x_1,x_2\bigr) \ldots p^{n+1,+}_{t_n-t_{n-1}}\bigl(x_{n-1},x_n\bigr).
\end{multline*}
\end{proof}

Notice that in the above proof, if we integrate $y_n$ over some smaller set than $W^n(x_n)$ we find that 
\begin{equation}
\label{filter}
{\mathbf Q}^{n,+}_{0,0} \bigl( Y_{t_n} \in A| X_{t_1},X_{t_2}, \ldots , X_{t_n} \bigr) = \int_{A \cap W^n(X_{t_n})} \lambda^n(X_{t_n}, 
y\bigr )dy.
\end{equation}
This may be interpreted as the following filtering property: the conditional distribution of $Y_t$ given $\bigl( X_s; s \leq t \bigr)$ is 
given by the density $\lambda^n(X_t, \cdot)$ on $W^n(X_t)$.

\section{Brownian motion in the Gelfand-Tsetlin cone}

Proposition \ref{interlace} lends itself to an iterative procedure. Let ${\mathbf K}$ be the cone of  points ${\mathbf x}= 
\bigl(x^{1},x^{2}, \ldots x^{n}\bigr)$ with  $x^{k}=\bigl(x^{k}_1,x^{k}_2, \ldots ,x^{k}_k\bigr) \in {\mathbf R}^k$ satisfying 
the inequalities 
\begin{equation}
 x^{k+1}_{i} \leq x^{k}_i \leq x^{k+1}_{i+1}.
\end{equation}
${\mathbf K}$ is sometimes called the Gelfand-Tsetlin cone, and arises in representation theory.
We will consider a process ${\mathbf X}(t)=\bigl( X^{1}(t), X^{2}(t), \ldots X^{n}(t)\bigr)$ taking values in ${\mathbf K}$ so that 
\begin{equation}
\label{gts:sde}
X^{k}_i(t)= x^{k}_i+\gamma^{k}_i(t)+ L^{k,-}_i(t)-L^{k,+}_i(t),
\end{equation}
where $\bigl(\gamma^k_i(t); t \geq 0 \bigr)$ for $1\leq k \leq n, 1 \leq i \leq k $ are independent Brownian motions, and 
$\bigl(L^{k,+}_i(t); t \geq 0 \bigr)$ and  $\bigl(L^{k,1}_i(t); t \geq 0 \bigr)$ are continuous, increasing processes growing only when 
$X^k_i(t)=X^{k-1}_i(t)$ and $X^k_i(t)=X^{k-1}_{i-1}(t)$ respectively,  the exceptional cases  $L^{k,+}_k(t)$ and $L^{k,-}_1(t)$  being 
identically zero for all $k$. 
For initial co-ordinates satisfying $x^{k}_i<x^{k}_{i+1}$ for all $k$ and $i$, we may give a pathwise construction, as in Section 2, based 
on alternately using the Skorokhod construction to reflect $X^{k}_i$ downwards from $X^{k-1}_i$ and upwards from $X^{k-1}_{i-1}$. The 
potential difficulty that $X^{k-1}_i$ meets  $X^{k-1}_{i-1}$ does not arise. 

In order to construct ${\mathbf X}$ starting from the origin we use a different method.
First we note that if the pair of processes $\bigl(X,Y\bigr)$, governed by the measure ${\mathbf Q}^{n,+}_{0,0}$ satisfies equations 
\eqref{eqndyson} and \eqref{eqnref}, then,
\begin{equation}
\label{obs}
\text{ the Brownian motion $\gamma= (\gamma_1, \gamma_2, \ldots \gamma_{n+1})$ is independent of $Y$.}
\end{equation}
By repeated application of Proposition \ref{interlace},  there exists a process $\bigl( {\mathbf X}(t); t \geq 0 \bigr)$, starting from the 
origin,  such that
\begin{description}
\item   the process $\bigl(X^{k}(t); t \geq 0 \bigr)$ is distributed as ${\mathbf P}^{k,+}_0$, for $k=1,2,\ldots, n$,
\item  the pair of processes $\bigl(X^{k+1}(t),X^{k}(t) ; t \geq 0 \bigr)$ are distributed as under ${\mathbf Q}^{k,+}_{0,0}$, for 
$k=1,\ldots, n-1$,
\item for $k=2, \ldots, n-1$ the process $\bigl(X^{k+1}(t); t \geq 0 \bigr)$  is conditionally independent of  $\bigl(X^1(t),\ldots, 
X^{k-1}(t); t \geq 0 \bigr)$ given $\bigl(X^{k}(t); t \geq 0 \bigr)$.
\end{description}
By its very construction the process ${\mathbf X}$  satisfies the equations \eqref{gts:sde}, for some  Brownian motions $\gamma^k_i$, which 
by the observation \eqref{obs} are independent.
Even starting from the origin, pathwise uniqueness, and hence uniqueness in law hold for ${\mathbf X}$.
Consequently we may state the  following proposition.

\begin{proposition}
\label{prop:gts}
The process $\bigl({\mathbf X}(t); t \geq 0 \bigr)$, satisfying \eqref{gts:sde},  if started from the origin, satisfies for each 
$k=1,2,\ldots ,n$, 
\[
\bigl( X^{(k)}(t); t \geq 0 \bigr) \text{ is distributed as under } {\mathbf P}^{k,+}_0.
\]
\end{proposition}

The conditional distribution of $\bigl(X^k(t); t \geq 0 \bigr)$ given  $\bigl(X^{k-1}(t); t \geq 0 \bigr)$ factorizes in a Markovian fashion 
into the product of the conditional distribution of $\bigl(X^k(s);  0\leq s \leq t \bigr)$ given  $\bigl(X^{k-1}(s);  0 \leq s \leq 
t\bigr)$, and the conditional distribution of $\bigl(X^k(u);  t\leq u  \bigr)$ given  $\bigl(X^{k-1}(u);  t \leq u \bigr)$ and $X^k(t)$. 
From this factorization  we deduce that, for any $t>0$, and $k=2, \ldots, n-1$ the process $\bigl(X^{k+1}(s); 0 \leq s \leq t \bigr)$ is 
conditionally independent of $\bigl(X^1(s),\ldots X^{k-1}(s); 0 \leq s \leq t \bigr)$ given $\bigl(X^{k}(s);  0 \leq s \leq t \bigr)$. For 
any $x^k\in W^k$ we will denote by ${\mathbf K}(x^k)$ the set of all $(x^1, x^2, \ldots x^{k-1})$  such that for all $i$ and $j$, 
$x^{j+1}_{i} \leq x^{j}_i \leq x^{j+1}_{i+1}$.
The $k(k-1)/2$-dimensional volume of ${\mathbf K}(x^k)$ is given by
\[
\frac{1}{\prod_{j<k} j!} h_{k}(x^k).
\]
  Recall from \eqref{filter} that the conditional distribution of 
$X^{k}(t)$ given $\bigl(X^{k+1}(s); 0 \leq s \leq t \bigr)$ has the density $\lambda^k(X^{k+1}(t), \cdot)$ on $W^{k}(X^{k+1}(t))$.  Combining 
this  with the conditional independence property we deduce that the  conditional distribution of $\bigl(X^1(t), X^2(t) \ldots 
X^{k}(t)\bigr)$  given $\bigl( X^{k+1}(s); 0 \leq s \leq t\bigr)$ is uniform on ${\mathbf K}(X^{k+1}(t))$. 
Finally using  the fact that the distribution of $X^n(t)$ is given by the density $\mu^n_t$  on $W^n$ we deduce that the distribution of  
${\mathbf X}(t)$ has the density
\begin{equation}
\label{entrance:gts}
\boldsymbol{\mu}^n_t({\mathbf x})= (2\pi)^{-n/2} t^{-n^2/2} \exp \left\{ -\sum_{i} (x^n_i)^2/(2t)\right\} \left\{ \prod_{i<j} (x^n_j-x^n_i) 
\right\},
\end{equation}
with respect to Lebesgue measure on ${\mathbf K}$.
Baryshnikov, \cite{baryshnikov}, studies this distribution  in some detail. Let  $\bigl( H(t); t \geq 0\bigr)$ be a Brownian motion in the 
space of $n \times n$ Hermitian matrices, and consider the process $\bigl(H^1(t), H^2(t), \ldots H^n(t); t \geq 0\bigr)$ where $H^k(t)$ is 
the $k\times k$ minor of $H(t)=H^n(t)$ obtain by deleting the last $n-k$ rows and columns. It is a classical result that the eigenvalues of 
$H^{k-1}(t)$ are interlaced with those of $H^k(t)$. Baryshnikov shows that, at any fixed instant $t>0$, the distribution of the   
eigenvalues of $H^1(t), H^2(t), \ldots H^n(t)$ is given by the density \eqref{entrance:gts}. However it is not the case that the eigenvalue 
process is distributed as the process $\bigl({\mathbf X}(t); t \geq 0 \bigr)$.

O'Connell, \cite{oconnell}, describes another process $\bigl(\boldsymbol{\Gamma}(t); t \geq 0 \bigr)$ taking values in ${\mathbf K}$ which 
is constructed via certain explicit  path transformations. This process arises as the scaling limit of the RSK correspondence. The process 
${\mathbf X}$ described above  has several features in common with $\boldsymbol{\Gamma}$. For each $k$, the subprocess $\bigl(\Gamma^k(t); 
t\geq 0\bigr)$ evolves as Dyson $k$-tuple starting from zero. Additionally 
\begin{equation}
\label{gammax}
\big( X^1_1(t), X^2_2(t),\ldots , X^n_n(t); t \geq 0 \bigr) \stackrel{dist}{=} \big( \Gamma^1_1(t), \Gamma^2_2(t),\ldots , \Gamma^n_n(t); t 
\geq 0 \bigr),
\end{equation}
but remarkably all other components $\Gamma^k_l$ with $l<k$ are given by explicit deterministic transformations applied to the processes 
$\Gamma^1_1, \Gamma^2_2, \ldots \Gamma^n_n$. A feature that ${\mathbf X}$ certainly does not share.

Notice that, for $k \geq 2$,
\begin{equation}
\label{xkk}
X^{k}_{k}(t)= \gamma_k^k(t)+ L^{k,-}_k(t),
\end{equation}
where $L^{k,-}_k(t)$ grows only when $X^k_k(t)=X^{k-1}_{k-1}(t)$. On applying the Skorokhod lemma, see Chapter VI of \cite{revuzyor}, we 
find that
\begin{equation}
L^{k,-}_k(t)= \sup_{s \leq t} \bigl( X^{k-1}_{k-1}(s)- \gamma^k_k(s)\bigr).
\end{equation}
Iterating this relation we obtain
\begin{equation}
X^k_k(t)= \sup_{0=t_0 \leq t_1 \leq t_2\leq \ldots \leq t_k=t} \sum_{i=1}^k \bigl\{\gamma^i_i(t_i)-\gamma^{i}_{i}(t_{i-1})\bigr\},
\end{equation}
which in the light of Proposition \ref{prop:gts} proves the identity \eqref{identity}.  This is  essentially the same argument for 
\eqref{identity} as given by O'Connell and Yor, \cite{oconnellyor},  with Proposition \ref{prop:gts} replacing the corresponding statement 
about $\boldsymbol{\Gamma}$. 
 
We close this section by noticing that   $\big( X^1_1(t), X^2_2(t),\ldots , X^n_n(t); t \geq 0 \bigr)$ is Markovian and giving an explicit 
formula for it transition probabilities.
For $n\geq 1$ let $\Phi^{(n)}_t$ denote the $n$th order iterated integral of the Gaussian density $\phi_t$,
\begin{equation}
\Phi^{(n)}_t(y)= \int_{-\infty}^y \frac{(y-x)^{n-1}}{(n-1)!} \phi_t(x)dx,
\end{equation}
and for $n \geq 0 $ let $\Phi^{(-n)}_t$ denote the $n$th order derivative of $\phi_t$.
Define for $x,x^\prime \in W^n$,
\begin{equation}
r_t(x,x^\prime)= \det \bigl\{ \Phi_t^{(i-j)}(x^\prime_j-x_i); 1 \leq i,j \leq n \bigr\}.
\end{equation}

\begin{lemma}
\label{weakcon2}
For any $f:W^{n} \rightarrow {\mathbf R}$ which is bounded and continuous and zero in a neighbourhood of the boundary of $W^n$,
\[
\lim_{t \downarrow 0} \int_{W^{n}} r_t\bigl( x,x^\prime\bigr)f(x^\prime)dx^\prime = f(x),
\]
 uniformly for all $x \in W^{n}$.
\end{lemma}

\begin{proposition}
The process $\big( X^1_1(t), X^2_2(t),\ldots , X^n_n(t); t \geq 0 \bigr)$ satisfying \eqref{xkk} is Markovian with transition densities 
given by $r_t(x,x^\prime).$
\end{proposition}
\begin{proof}
For a fixed $x^\prime \in {\mathbf R}$, and any $n$, the function $(t,x) \mapsto \Phi^{(n)}_t(x^\prime-x)$ solves the heat equation on 
$(0,\infty) \times {\mathbf R}$. From this we easily see that for a fixed $x^\prime \in {\mathbf R}^n$, that the function $(t,x) \mapsto 
r_t(x^\prime-x)$ solves the heat equation on $(0,\infty) \times {\mathbf R}^n$. Moreover if $x_i=x_{i-1}$ for any $i=2,3, \ldots, n$ then 
the $i$th and $(i-1)$th rows of the determinant defining $\frac{\partial}{\partial x_i}r_t(x,x^\prime)$ are equal and hence this quantity is 
zero.

Let $f: W^{n} \rightarrow {\mathbf R}$ be a bounded, continuous and are in a neighbourhood of the boundary of $W^n$. Then   define a smooth 
function $F$  on $(0,\infty)\times W^{n}$ via
\begin{equation*}
F(t,x)= \int_{W^{n}} r_t (x, x^\prime) f(x^\prime)dx^\prime.
\end{equation*}
By virtue of the above observations regarding $r_t$, and differentiating through the integral,  we find that
\begin{equation*}
\tfrac{1}{2} \sum_{i=1}^{n} \frac{\partial^2 F}{\partial x_i^2} (t,x)
 = \frac{\partial F}{\partial t}(t, x)   \quad \text{ on } (0,\infty)\times W^{n},
\end{equation*}
with the boundary conditions
\begin{equation*}  
\frac{\partial F}{\partial x_i}(t,x)=0 \text{  whenever  
$x_i=x_{i-1}$ for some $i=2,3,\ldots ,n$.} 
\end{equation*}
Let $X$ denote a process governed by a probability ${\mathbf R}_x$, with components $X_1(t) \leq  X_2(t)\leq  \ldots \leq X_n(t)$ satisfying 
the equations
$X_k(t)=x_k +\gamma_k(t) + L^k(t)$, where $\gamma_k$ are independent Brownian motions and $L^k$ is an increasing process growing only when 
$X_k(t)=X_{k-1}(t)$, with $L^1$ being identically zero.
Fix $T, \epsilon>0$. Applying It\^{o}'s formula, we find that the process $\bigl(F\bigl(T+\epsilon-t, X_t\bigr); t\in [0,T]\bigr)$ is a 
local martingale, which being bounded is a true martingale. Thus 
\begin{equation*}
F\bigl(T+\epsilon, (x) \bigr)={\mathbf R}_{x} \bigl[ F\bigl(\epsilon,  X(T)\bigr) \bigr].
\end{equation*}
Appealing to the previous lemma, we may let $\epsilon \downarrow 0$ and so obtain,
\[
F(T, x)= {\mathbf R}_{x} \bigl[ f( X(T)) \bigr],
\]
which, since it is clear the distribution of $X(T)$ does not charge the boundary of $W^n$,  proves the proposition.
\end{proof}

In view of Proposition \ref{prop:gts}, we obtain from $r_t$ by a simple integration the following expression for the distribution function 
of the largest eigenvalue of $H(t)$:
\begin{equation}
\label{topdist}
{\mathbf P}^{n,+}_{0}\bigl( X_n \leq x \bigr)= \det\bigl\{ \Phi_t^{(i-j+1)}(x)\bigr\}.
\end{equation}
Possibly this can be checked directly using the Heine identity.

\section{Coalescing Brownian motions}

In this section we consider  the joint distribution of a family of coalescing Brownian motions. 
Fix  $z_1 \leq z_2\leq \ldots \leq z_n$  and consider the process of $n$ coalescing Brownian Motions, 
\[
t \mapsto Z_t= \bigl( Z_t(z_1), \ldots Z_t(z_n)\bigr),
\]
where each process $\bigl( Z_t(z_i); t \geq 0 \bigr)$ is a Brownian motion (relative to some common filtration) starting from 
$Z_0(z_i)=z_i$, with for each 
distinct pair $i\neq j$  the process
\[
t \mapsto \frac{1}{\sqrt{2}} | Z_t(z_i)-Z_t(z_j)|
\]
being a standard Brownian motion on the half-line $[0,\infty)$ with an absorbing barrier at $0$.
Thus informally $\bigl(Z_t(z_i); t \geq 0 \bigr)$ and $\bigl(Z_t(z_j); t\geq 0 \bigr)$ evolve independently until they first meet, after 
which they coalesce and move together. Such families of coalescing Brownian motions have been well-studied, for some recent works concerning 
them see 
\cite{fontesetal} and \cite{evanszhou}.

For a fixed $t>0$, the distribution of $Z_t(z)$ is supported on $W^n$. That part of the distribution supported on the boundary of $W^n$ 
corresponds to the event that coalescence has occurred. Whereas the restriction of the distribution to the interior $W^n$  (corresponding  
to no coalescence) is  given by Karlin-McGregor formula :
\begin{equation}
{\mathbf P}\bigl( X_t(z_i) \in dz^\prime_i \text{ for all $i$}\bigr) = \det \bigl\{ \phi_t(z^\prime_j-z_i)\bigr\} dz^\prime.
\end{equation}
In fact we can bootstrap from  this  result to a complete determination of the law of $Z_t(z)$, which can be expressed in the following neat 
way. 

\begin{proposition}
\label{distribution}
For $z,z^\prime \in {W}^n$, the probability
\[
 {\mathbf P}\bigl( Z_t(z_i)\leq z^\prime_i \text{ for $1\leq i\leq n$}\bigr)\]
is given by the determinant of an $n\times n$ matrix  with $(i,j)$th element given by 
\begin{align*}
\Phi_t(z^\prime_j-z_i) & \qquad \text{ if } i \geq j, \\
\Phi_t(z^\prime_j-z_i) -1 & \qquad  \text{ if } i<j,
\end{align*}
where
\[
\Phi_t(z)=\int_{-\infty}^z \frac{dy}{\sqrt{2\pi t}} \exp\{ -y^2/(2t)\}.\]
\end{proposition}
\begin{proof} 
First we note that by integrating the Karlin-McGregor formula we obtain
\begin{equation}
{\mathbf P}\bigl( Z_t(z_1) \leq z^\prime_1 < X_t(z_2) \leq z^\prime_2< \ldots 
\leq z^\prime_{n-1}<Z_t(z_n) \leq z^\prime_n\bigr)= \det\bigl\{ \Phi_t(z^\prime_j-z_i)\bigr\}.
\end{equation}
We are going to obtain the desired result by showing how the indicator function of the event of interest
\[
\bigl\{ Z_t(z_1) \leq z^\prime_1, Z_t(z_2) \leq z^\prime_2, \ldots, Z_t(z_n) \leq z^\prime_n \bigr\}
\]
can be expanded in terms of the indicator functions of the events of the form
\[
\{ Z_t(z_{i(1)}) \leq z^\prime_{j(1)} < Z_t(z_{i(2)}) \leq z^\prime_{j(2)} < \ldots < 
z^\prime_{j(s-1)}<Z_t(z_{i(s)}) \leq z^\prime_{j(s)}\bigr\},
\]
for increasing subsequences of indices $i(1),i(2), \ldots ,i(s)$ and $j(1), j(2), \ldots ,j(s)$.
To this end I claim  firstly that, whenever $z, z^\prime\in {W}^n$,
\begin{equation}
\det\bigl\{ {\bf 1}(z_i\leq  z^\prime_j )\bigr\}= {\bf 1} (z_1\leq z^\prime_1 < z_2 \leq z^\prime_2<\ldots < z_n \leq z^\prime_n). 
\end{equation}
I claim secondly that
\begin{equation}
\det \begin{Bmatrix}
{\bf 1} (z_i \leq z^\prime_j) &  i \geq j \\
-{\bf 1}(z^\prime_j<z_i) & i <j
\end{Bmatrix}
={\bf 1}( z_1 \leq  z^\prime_1, z_2\leq z^\prime_2, \ldots ,z_n\leq z^\prime_n).
\end{equation}

To prove the first claim take the matrix $M=\bigl\{ {\bf 1}(z_i\leq  z^\prime_j )\bigr\}$, and subtract from each column (other than the 
first) the values of the preceding column.  The diagonal elements  of this new matrix are 
\[
{\mathbf 1}( z_i \leq z^\prime_i)- {\mathbf 1}(z_i\leq z^{\prime}_{i-1})=
{\mathbf 1}( z_{i-1}^\prime< z_i \leq z_i^\prime);
\]
adopting the convention that $z^\prime_0=-\infty$. Thus the product of these diagonal elements gives the desired result. We have to check 
that in the expansion of the determinant this is the only contribution. Suppose that $\rho$ is a permutation, not the identity. Then we can 
find $i<j$ with $\rho(i)>i$ and $\rho(j)\leq i$. Consider the product of the $(i,\rho(i))$th and $(j,\rho(j))$th elements of the matrix 
(after the column operations). We obtain
\[
{\mathbf 1} ( z^\prime_{\rho(i)-1} <z_i \leq z^\prime_{\rho(i)} ) {\mathbf 1}
 ( z^\prime_{\rho(j)-1} <z_j \leq z^\prime_{\rho(j)} ).
\]
This can only be non-zero if both $z^\prime_{\rho(i)-1}< z_i$ and $z_j \leq z^\prime_{\rho(j)} $; but $z_i\leq z_j$ so this would imply 
$z^\prime_{\rho(i)-1}<
z^\prime_{\rho(j)}$. In view of the fact  $\rho(i)-1 \geq \rho(j)$ this is impossible. 

Consider the matrix $N$  appearing in the second claim. The product of its diagonal elements gives the desired result. To show that this is 
the only contribution to the determinant, take $\rho$ a permutation, not equal to the identity and $i<j$ with $\rho(i)>i$ and $\rho(j) \leq 
i$, as before. Then the product of the $(i, \rho(i))$th and $(j,\rho(j))$th elements of the matrix is 
\[
-{\mathbf 1}(z^\prime_{\rho(i)}<z_i) {\mathbf 1}( z_j \leq z^\prime_{\rho(j)})
\]
Since $z_i\leq z_j$ for this to be non-zero we would have to have $ z^\prime_{\rho(i)} < z^\prime_{\rho(j)}$, which is impossible for 
$\rho(i)> \rho(j)$. 

Let $T=\{ -{\mathbf 1}(j>i)\}$ be the upper triangular matrix so that $N=M+T$ and consider the Laplace expansion of $\det(M+T)$ in terms of 
minors. For  increasing  vectors of subscripts ${\mathbf i}$ and ${\mathbf j}$  let $M[{\mathbf i},{\mathbf j}]$ denote the corresponding 
minor of $M$ and let $\tilde{T}[{\mathbf i},{\mathbf j}]$ be the complementary minor of $T$ so that
\[
\det(N)= \det(M+T)= \sum_{{\mathbf i},{\mathbf j}} (-1)^{s({\mathbf i},{\mathbf j})} M[{\mathbf i},{\mathbf j}]\tilde{T}[{\mathbf 
i},{\mathbf j}],
\]
for appropriate signs $s({\mathbf i},{\mathbf j})$.
Evaluating $\det(N)$ via the second claim, and the minors $ M[{\mathbf i},{\mathbf j}]$ via (general versions of ) the first claim we have 
obtained an expansion of ${\bf 1}( z_1 \leq  z^\prime_1, z_2\leq z^\prime_2, \ldots ,z_n\leq z^\prime_n)$ as a linear combination of of 
terms of the form ${\mathbf 1}( z_{i(1)} \leq z^\prime_{j(1)} < z_{i(2)} \leq z^\prime_{j(2)} < \ldots \leq z^\prime_{j(s)})$. 

To complete the proof replace, in the above expansion, $z_i$ by $Z_t(z_i)$ and take expectations. On the lefthandside we obtain ${\mathbf 
P}\bigl( Z_t(z_1) \leq z^\prime_1, Z_t(z_2) \leq z^\prime_2, \ldots, Z_t(z_n) \leq z^\prime_n \bigr)$. On the righthandside we  have a 
linear combination of probabilities:
${\mathbf P}\bigl( Z_t(z_{i(1)}) \leq z^\prime_{j(1)} < Z_t(z_{i(2)}) \leq z^\prime_{j(2)} < \ldots < 
z^\prime_{j(s-1)}<Z_t(z_{i(s)}) \leq z^\prime_{j(s)}\bigr)$ each of which can re-written by means of the integrated Karlin-McGregor formula 
as a minor of the 
determinant $\det\{ \Phi(z^\prime_j-z_i) \}$. And to finish we notice that the righthandside is now the Laplace expansion of the   
determinant of the sum of matrices $\{ \Phi(z^\prime_j-z_i)\}$ and $\{ - {\mathbf 1}(j>i)\}$.
\end{proof}

The expression just obtained  for the distribution of coalescing Brownian motions is closely related to the formula  for the transition 
density 
of the interlaced Brownian motions given  by Proposition \ref{trans1}. In fact it is easily verified that
\begin{multline}
\label{dual1}
q^n_t\bigl((x,y), (x^\prime,y^\prime) \bigr)= 
(-1)^n\frac{\partial^n}{\partial y_1\ldots \partial y_n} \frac{\partial^{n+1}}{\partial x^\prime_1\ldots \partial x^\prime_{n+1} }{\mathbf 
P} \bigl( Z_t(x_i) \leq x^\prime_i, Z_t(y_j) \leq y^\prime_j \text{ for all } i,j \bigr).
\end{multline}
This represents a duality between the the interlaced  Brownian motions and coalescing  Brownian motions which generalizes the well-known 
duality between Brownian motion on the half-line $[0,\infty)$ with a reflecting Barrier at zero, and Brownian motion on the half-line with 
an absorbing barrier at zero.

There is interesting alternative  way of expressing the equality \eqref{dual1}. The Arratia flow or Brownian web   is a infinite family of 
coalescing Brownian motions,
with a path starting from every point in space-time. Let $ t\in [s,\infty) \mapsto Z_{s,t}(x)$ denote the path starting from $(s,x)$. It is 
possible to define on  the same  probability space a dual flow with paths running backwards in time: $ s \in (-\infty,t] \mapsto 
\hat{Z}_{s,t}(x)$ being the path beginning at $(t,x)$. For the details of this construction see \cite{tothwerner} and \cite{fontesetal}. The 
 flow $Z$ and its dual $\hat{Z}$ are such that for any $s,t,x$ and $y$,  the two events $ Z_{s,t}(x) \leq y$ and $\hat{Z}_{s,t}(y) \geq x $ 
differ by a set of zero probability. Using this we may rewrite \eqref{dual1} as
\begin{equation}
\label{dual2}
q^n_t\bigl((x,y), (x^\prime,y^\prime) \bigr)dx^\prime dy = {\mathbf P} \bigl( Z_{0,t}(x_i) \in dx^\prime_i, \hat{Z}_{0,t}(y^\prime_j) \in 
dy_j \text{ for all } i,j \bigr).
\end{equation}
This seems to fit with the  fact that the paths of $\hat{Z}$ are ''reflected off'' those of $Z$, see \cite{soucaliuctothwerner} and 
\cite{sw}.  

\section{Proofs of two lemmas}

\begin{proof}[Proof of Lemma \ref{weakcon2}]
The contribution to the determinant defining $r_t(x,x^\prime)$ coming from the principal diagonal is equal to the standard heat kernel in 
${\mathbf R}^n$. The lemma will follow if we can show all other contributions to the determinant are uniformly negligible as $t$ tends down 
to $0$. 
Choose $\epsilon>0$ so that the function $f$ is zero in an $2\epsilon$-neighbourhood of the boundary of $W^n$.  Then   consider a 
contribution to the determinant corresponding to some permutation $\rho$ which is not the identity. There exist  $i<j$ with $\rho(i)>i$ and 
$\rho(j)\leq i$, and the contribution corresponding to $\rho$ consequently contains factors of 
$\Phi_t^{(i-\rho(i))}(x^\prime_{\rho(i)}-x_i)$ and $\Phi_t^{(j-\rho(j))}(x^\prime_{\rho(j)}-x_j)$.  Noting that $j-\rho(j)>0$ and 
$i-\rho(i)<0$ we see that on the set $\{x^\prime_{\rho(i)}-x_i> \epsilon\} \cup \{x^\prime_{\rho(j)}-x_j<-\epsilon\}$ at least one of these 
factors, and indeed the entire contribution, tends to zero uniformly  as $t$ tends down to zero.  But on the complement of this set we have 
$x^\prime_{\rho(i)}\leq x_i+\epsilon \leq x_j+\epsilon \leq x^\prime_{\rho(j)}+2\epsilon$,
and $\rho(j)\leq \rho(i)$ implies that $x^\prime_{\rho(j)} \leq x^\prime_{\rho(i)}$, so we see that $x^\prime$ is within the 
$2\epsilon$-neighbourhood  of the boundary of $W^n$, and does not belong to the support of $f$. This proves the lemma.
\end{proof}

\begin{proof}[Proof of Lemma \ref{weakcon}]
It is convenient to write $z_1=x_1, z_3=x_2, \ldots z_{2n+1}=x_n$, and  $z_2=y_1, z_4=y_2, \ldots, z_{2n}=y_n$, with a corresponding change 
of notation for $x^\prime_i$ and $y^\prime_i$ also. Now  reorder the columns and rows of the determinant defining $q^n_t$ so that the 
$(i,j)$th entry is a function of the difference $z^\prime_j-z_i$. We may now argue in the same way as in the preceding proof. Choose 
$\epsilon>0$ so that the function $f$ is zero in an $2\epsilon$-neighbourhood of the boundary of $W^{n+1,n}$. Consider a contribution to the 
determinant corresponding to some permutation $\rho$ which is not the identity. There exist  $i<j$ with $\rho(i)>i$ and $\rho(j)\leq i$, and 
the contribution corresponding to $\rho$ consequently contains factors  which are functions of $z^\prime_{\rho(i)}-z_i$ and 
$z^\prime_{\rho(j)}-z_j$.  Noting that $j-\rho(j)>0$ and $i-\rho(i)<0$, and checking the entries of the determinant above and below the 
diagonal we see that on the set $\{z^\prime_{\rho(i)}-z_i> \epsilon\} \cup \{z^\prime_{\rho(j)}-z_j<-\epsilon\}$ at least one of these 
factors, and indeed the entire contribution, tends to zero uniformly  as $t$ tends down to zero. As above, this proves the lemma.
\end{proof}

\vspace{.1in}

Department of Statistics, University of Warwick, Coventry CV4 7AL, UK.
\vspace{.1in}

{\em E-mail address:}  {\tt j.warren@warwick.ac.uk}

\end{document}